\documentclass[12pt,a4paper]{article} 
\usepackage{a4}

\usepackage[intlimits,tbtags]{amsmath} 
\usepackage{amsfonts,amssymb,amsthm}

\usepackage{cite,latexsym}

\newcommand{\tr}{\triangleright}
\newcommand{\tl}{\triangleleft} 
\newcommand{\R}{\mathcal{R}}
\newcommand{\RI}{{\mathcal{R}_{\mathrm{I}}}} 
\newcommand{\RII}{\mathcal{R}_{\mathrm{II}} } 

\newcommand{\op}{\mathrm{op}} 
\newcommand{\adL}{{\mathrm{ad_L}}}
\newcommand{\adR}{{\mathrm{ad_R}}}

\newcommand{\Acal}{\mathcal{A}} 
\newcommand{\Hcal}{\mathcal{H}}
\newcommand{\Xcal}{\mathcal{X}} 
\newcommand{\Lcal}{\mathcal{L}}

\newcommand{\Ket}[1]{\lvert #1\rangle}

\newcommand{\lrAngle}[1]{\langle #1\rangle}

\newtheorem{Definition}{Definition} 
\newtheorem{Theorem}{Theorem}
\newtheorem{Proposition}{Proposition}

\newcommand{\uq}{{\mathcal{U}_q(\mathrm{u}_1)}}

\newcommand{\su}{{\mathcal{U}(\mathrm{su}_2)}}

\newcommand{\suq}{{\mathcal{U}_q(\mathrm{su}_2)}}
\newcommand{\slC}{{\mathcal{U}_q(\mathrm{sl}_2(\mathbb{C})) }}

\newcommand{\SUq}{{SU_q(2)}} 

\newcommand{\Mink}{{\mathbb{R}_q^{1,3}}}

\newcommand{\qsphere}{{\mathcal{S}_{q\infty}^\op}}

\newenvironment{Proplist}[1]%
{\begin{list}{}{ \settowidth{\labelwidth}{\textrm{#1}}
      \setlength{\leftmargin}{\parindent}
      \addtolength{\leftmargin}{\labelwidth}
      \addtolength{\leftmargin}{\labelsep}}
      }%
    {\end{list}}

\begin{document}

\rightline{LMU-TPW 2001-12}

\vspace{4em}
\begin{center}

{\Large{\bf Spin in the \textit{q}-Deformed Poincar\'e Algebra}}

\vspace{3em}

\textbf{Christian Blohmann}

\vspace{1em}
 
Ludwig-Maximilians-Universit\"at M\"unchen, Sektion Physik\\
Lehrstuhl Prof.\ Wess, Theresienstr.\ 37, D-80333 M\"unchen\\[1em]

Max-Planck-Institut f\"ur Physik, 
        F\"ohringer Ring 6, D-80805 M\"unchen\\[1em]

\end{center}

\vspace{1em}

\begin{abstract}
  We investigate spin as algebraic structure within the $q$-deformed
  Poin\-ca\-r\'e algebra, proceeding in the same manner as in the
  undeformed case.  The $q$-Pauli-Lubanski vector, the $q$-spin
  Casimir, and the $q$-little algebras for the massless and the
  massive case are constructed explicitly.
\end{abstract}

\section{Introduction}

From the beginnings of quantum field theory it has been argued that
the pathological ultraviolet divergences should be remedied by
limiting the precision of position measurements by a fundamental
length \cite{Born:1933,March:1936,Heisenberg:1938}. In view of how
position-momentum uncertainty enters into quantum mechanics, a natural
way to integrate such a position uncertainty in quantum theory would
have been to replace the commutative algebra of space observables with
a non-commutative one \cite{Snyder:1947}. However, deforming the space
alone will in general break the symmetry of spacetime. In order to
preserve a background symmetry the symmetry group must be deformed
together with the space it acts on. This reasoning led to the
discovery of quantum groups \cite{Drinfeld:1986}, that is, generic
methods to continuously deform Lie algebras
\cite{Drinfeld:1985,Jimbo:1985} and matrix groups
\cite{Woronowicz:1987,Faddeev:1990,Takeuchi:1990} within the category
of Hopf algebras. Starting from the non-commutative plane
\cite{Manin:1988}, the $q$-deformations of a series of objects,
differential calculi on non-commutative spaces \cite{Wess:1991},
Euclidean space \cite{Faddeev:1990}, Minkowski space
\cite{Carow-Watamura:1990}, the Lorentz group and the Lorentz algebra
\cite{Podles:1990,Carow-Watamura:1991,Schmidke:1991,Ogievetskii:1991a},
led to the $q$-de\-formed Poincar\'e algebra
\cite{Ogievetskii:1992a,Majid:1993}.

Describing the symmetry of flat spacetime, the Poin\-ca\-r\'e group
or, equivalently, its enveloping algebra is sufficient to construct
special relativity and even a considerable part of relativistic
quantum theory. More precisely, Wigner has shown that a free
elementary particle can be identified with an irreducible Hilbert
space representation of the Poincar\'e group \cite{Wigner:1939}. These
representations are constructed using the method of induced
representations, which reduces the representation theoretic problem to
a structural analysis of the Poincar\'e algebra. In mathematical
terms: We start from a representation of the inhomogeneous part of the
algebra, determine the stabilizer (little algebra) of this
representation, construct the irreducible representations of the
stabilizer, and, finally, induce these representations to
representations of the entire algebra, yielding \emph{all} irreducible
representations of the Poincar\'e algebra.

Seemingly abstract, each step in this construction has a clear
physical interpretation: The representation of the inhomogeneous part
is the description of a momentum eigenspace. The stabilizer is the
spin symmetry lifting a possible momentum degeneracy. The
representation of the stabilizer defines the transformations of the
spin degrees of freedom, the canonical example being a massive spin
particle at rest carrying a representation of $\mathrm{SU}(2)$.
Finally, the induction is the boosting of a rest state to arbitrary
momentum. We see that this procedure is not only the mathematical
means to construct the wanted representations, but provides insight in
the physical nature of spin. It tells us that there is spin, because
in general momentum is not sufficient to characterize a particle
uniquely. It tells us what the symmetry structure of the spin degrees
of freedom is, $\mathrm{SU}(2)$ in the massive case but
$\mathrm{ISO}(2)$ in the massless case. And it tells us, that momentum
and spin are \emph{all} possible exterior degrees of freedom of a
particle.

The physical line of thought described in the last paragraph relies on
the sole assumption that the Poincar\'e algebra describes the basic
symmetry of spacetime. The $q$-deformed Poincar\'e algebra has been
constructed to describe the basic symmetry of $q$-deformed space time.
Therefore, we can proceed in exactly the same manner to find out what
$q$-deformed spin is.

In Sec.~\ref{sec:qPoin} we review the $q$-Poincar\'e algebra with
focus on its general structure. In Sec.~\ref{sec:qPL} we define the
key properties of a useful $q$-deformed Pauli-Lubanski vector and
present such a vector in theorem~\ref{th:PL}. Its square yields the
spin Casimir. Sec.~\ref{sec:qLittle} uses this $q$-Pauli-Lubanski
vector to compute the $q$-little algebras for both the massive and the
massless case.

Throughout this article, it is assumed that $q$ is a real number
$q>1$. We will frequently use the abbreviations $\lambda = q - q^{-1}$
and $[j] = \frac{q^j - q^{-j}}{q-q^{-1}}$ for a real number $j$, in
particular $[2] = q+q^{-1}$. The lower case Greek letters $\mu$,
$\nu$, $\sigma$, $\tau$ denote 4-vector indices running through
$\{0,-,+,3\}$. The upper case Roman letters $A$, $B$, $C$ denote
3-vector indices running through $\{-1,0,+1\} = \{-,3,+\}$.

\section{The \textit{q}-Deformed Poincar\'e Algebra}
\label{sec:qPoin}

The $q$-Poincar\'e algebra can be defined very explicitly by listing
its generators and the commutation relations between them. This has
been done in Appendix~\ref{sec:AppPoin}. Here, we give an overview of
the more general algebraic structure.

The \textbf{\textit{q}-Lorentz algebra} $\Hcal = \slC$ is a Hopf-$*$
algebra, with coproduct $\Delta$, counit $\varepsilon$, and antipode
$S$. We will also use the Sweedler notation $\Delta(h) =
h_{(1)}\otimes h_{(2)}$. Several forms of the $q$-Lorentz algebra can
be found in the literature, which are essentially equivalent. Here, it
is natural to use the form, where $\Hcal$ is described as Drinfeld
double of $\suq$ with its dual $\SUq^\op$ \cite{Podles:1990},
\begin{equation}
  \Hcal = \suq \Join \SUq^\op \,,
\end{equation}
that is, the Hopf-$*$ algebra generated by the algebra of rotations
$\suq$ and the algebra of boosts $\SUq^\op$ with cross commutation
relations
\begin{equation}
  bl = \lrAngle{ l_{(1)}, b_{(1)} } \,l_{(2)} b_{(2)} \,
  \lrAngle{S(l_{(3)}), b_{(3)}}
\end{equation}
for all $l\in\suq$, $b\in\SUq^\op$, where $\lrAngle{l,b}$ denotes the
dual pairing. In addition to the Drinfeld-Jimbo generators $E$, $F$,
$K=q^H$ of $\suq$ we will also use the the Casimir operator $W$ and
the 3-vector $\{J_A\} = \{J_-,J_+,J_3\}$ of angular momentum. The
generators $a$, $b$, $c$, $d$ of boosts form a multiplicative quantum
matrix $(\begin{smallmatrix} a&b\\c&d\end{smallmatrix})$. $\Hcal$
possesses two universal $\R$-matrices, $\RI$ and $\RII$, the first of
which is antireal $\RI^{*\otimes *} = \RI^{-1}$, the second is real
$\RII^{*\otimes *} = \R_{\mathrm{II}\,21}$. We often write in a
Sweedler like notation $\R = \R_{[1]} \otimes \R_{[2]}$.

The \textbf{\textit{q}-Minkowski space algebra} $\Xcal = \Mink$ is
generated by the 4-momentum vector $\{P_\mu\} =\{P_0, P_-, P_+,
P_3\}$ with relations
\begin{equation}
\label{eq:MinkRel}
    P_{\mu}P_{\nu} R_\mathrm{I}^{\nu\mu}{\!}_{\sigma\tau}
    = P_\sigma P_\tau \Leftrightarrow P_{\mu}P_{\nu}
    (R_\mathrm{I}^{-1})^{\mu\nu}{\!}_{\sigma\tau} = P_\tau P_\sigma \,,
\end{equation}
where the $R$-matrix $R_\mathrm{I}^{\mu\nu}{\!}_{\sigma\tau} =
\Lambda(\R_{\mathrm{I}[1]})^{\mu}{}_\sigma
\Lambda(\R_{\mathrm{I}[2]})^{\nu}{}_{\tau}$ is the 4-vector
representation of $\RI$. The 4-mo\-men\-tum vector is the basis of
this 4-vector representation of $\Hcal$, $h\tr P_\nu \equiv P_\mu
\Lambda(h)^\mu{}_\nu$, where $\Lambda$ is the representation map.
Relations~\eqref{eq:MinkRel} are the only homogeneous commutation
relations of $\Xcal$, which are consistent with this representation
and which have the right commutative limit. Consistency means that
$\Xcal$ is a left $\Hcal$-module $*$-algebra, that is,
\begin{xalignat}{2}
  h\tr xx' &= (h_{(1)}\tr x)(h_{(2)}\tr x')\,,
  & (h \tr x)^* &= (Sh)^* \tr x^*
\end{xalignat}
for all $h\in\Hcal$, $x\in\Xcal$.

The \textbf{\textit{q}-Poincar\'e algebra} $\Acal$ is the Hopf
semidirect product
\begin{equation}
  \Acal = \Xcal\rtimes\Hcal \,,
\end{equation}
the $*$-algebra generated by the $*$-algebras $\Xcal$ and $\Hcal$ with
cross commutation relations $hx = (h_{(1)}\tr x) h_{(2)}$. More
accurately, we have the following
\begin{Definition}
  Let $\Hcal$ be a Hopf-$*$ algebra and $\Xcal$ a left $\Hcal$-module
  $*$-algebra. The semidirect product $\Xcal\rtimes\Hcal$ is the
  $*$-algebra defined as the vector space $\Xcal\otimes\Hcal$ with
  multiplication
  \begin{equation}
  \label{eq:PoincCommute}
    (x\otimes h)(x'\otimes h') := x(h_{(1)}\tr x')\otimes h_{(2)}h'
  \end{equation}
  and $*$-structure $(x\otimes h)^* = (1\otimes h^*)(x^*\otimes 1)$.
  We often abbreviate $x \equiv x\otimes 1$ and $h \equiv 1\otimes h$.
\end{Definition}
There is a left and a right Hopf adjoint action of $\Hcal$ on $\Acal$
defined as
\begin{xalignat}{2}
  \label{eq:HopfAdjointAction}
  \adL h\tr a &:= h_{(1)} a S(h_{(2)}) \,,& a \tl \adR h &:=
  S(h_{(1)}) a h_{(2)} \,.
\end{xalignat}
The commutation relations~\eqref{eq:PoincCommute} are precisely such
that the left Hopf adjoint action of $\Hcal$ on $\Xcal$ equals the
module action $\adL h \tr x = h\tr x$. Let $\rho$ be a finite
representation of the $q$-Lorentz algebra $\Hcal$. We call a set of
operators $\{T_i\}$ a left or a right $\rho$-tensor operator if
\begin{equation}
\label{eq:TensorOp}
  \adL h\tr T_j = T_i\, \rho(h)^{i}{}_j \quad\mathrm{or}\quad
  T_j \tl \adL h = T_i\, \rho(S^{-1}h)^{i}{}_j
\end{equation}
holds, respectively, for all $h\in\Hcal$. By definition of the
$q$-Poincar\'e algebra, the momenta $P_\mu$ form a left
$\Lambda$-tensor operator, that is, a left 4-vector operator.

\section{The \textit{q}-Pauli-Lubanski Vector and the Spin Ca\-si\-mir}
\label{sec:qPL}

\subsection{Defining Properties of the  Pauli-Lubanski Vector}

In the undeformed case one defines the Pauli-Lubanski (pseudo)
4-vector operator
\begin{equation}
\label{eq:PauliLubanski1}
  W_\mu^{q=1} := -\frac{1}{2}\, \varepsilon_{\mu\nu\sigma\tau}
  L^{\nu\sigma} P^\tau \,,
\end{equation}
where $\varepsilon$ is the totally antisymmetric tensor,
$L^{\nu\sigma}$ the matrix of Lorentz generators, and $P^\tau$ the
momentum 4-vector. It is useful because each component of
$W_\mu^{q=1}$ commutes with each component of $P^\nu$, from which
follows that the 4-vector square $W^2 = \eta^{\mu\nu}
W_\mu^{q=1}W^{q=1}_\nu$ is a Casimir operator. The eigenvalues of this
Casimir operator are $-m^2s(s+1)$ where $s$ is the spin. Therefore,
$W_\mu^{q=1}$ can be viewed as square root of the spin Casimir.

In the $q$-deformed case we can try to define $W_\mu$ by
Eq.~\eqref{eq:PauliLubanski1}, as well, with the $q$-deformed versions
of the epsilon tensor, of the matrix of Lorentz generators, and of the
momenta. By construction, this definition yields a left 4-vector
operator. But the square of this 4-vector does not commute with the
momenta and, hence, it is not the searched-for spin Casimir.

In general, the assumption that $W_\mu$ and $P_\nu$ commute is not
consistent with both, $W_\mu$ and $P_\nu$, being \emph{left} 4-vector
operators. Otherwise, the expression $\adL h \tr [W_\mu,P_\nu]$ would
have to vanish for all $h\in\Hcal$, that is,
\begin{equation}
\label{eq:WPdep}
  W_{\mu'}P_{\nu'} \bigl( \Lambda(h_{(1)})^{\mu'}{\!}_\mu
  \Lambda(h_{(2)})^{\nu'}{\!}_\nu - \Lambda(h_{(2)})^{\mu'}{\!}_\mu
  \Lambda(h_{(1)})^{\nu'}{\!}_\nu \bigr)
  \stackrel{!}{=} 0 \,.
\end{equation}
Since the coproduct is not cocommutative as in the undeformed case,
this seems only possible for degenerate forms of $W_\mu$. We can avoid
this problem if we assume that the $q$-Pauli-Lubanski vector is a
\emph{right} 4-vector operator, making the following general
observation:
\begin{Proposition}
  \label{th:PL1}
  Let $a\in \Acal=\Xcal\rtimes \Hcal$ commute with $\Xcal$, $[a,x]=0$
  for all $x\in\Xcal$. Then $a\tl\adR h$ also commutes with $\Xcal$
  for any $h\in\Hcal$. In other words, the centralizer of $\Xcal$ is
  invariant under the right Hopf adjoint action of $\Hcal$.
\end{Proposition}
\begin{proof}
  Let $x\in\Xcal$ be any element of the quantum space. Then
\begin{equation}
\begin{split}
  (a \tl\adR h)\, x
  &= S(h_{(1)})a h_{(2)} x = S(h_{(1)})a (h_{(2)}\tr x) h_{(3)} \\
  &= S(h_{(1)}) (h_{(2)}\tr x) a h_{(3)} \\
  &= (S(h_{(1)})_{(1)} h_{(2)}\tr x) S(h_{(1)})_{(2)} a h_{(3)} \\
  &= (S(h_{(2)}) h_{(3)}\tr x) S(h_{(1)}) a h_{(4)} \\
  &= x \,(a \tl\adR h)
\end{split}
\end{equation}
for any $h\in\Hcal$.
\end{proof}
This means, that if a single component of a \emph{right} vector
operator commutes with all momenta, then the other components commute
with all momenta, as well. Hence, the requirement that $P_\mu$ and
$W_\nu$ commute does no longer generate linear dependencies of
type~\eqref{eq:WPdep}. We come to the following
\begin{Definition}
  \label{th:PLdef}
  A set of operators $W_\mu\in\Acal$ with the properties
  \begin{Proplist}{(PL3)}
  \item[\textup{(PL1)}] $W_\mu$ is a \emph{right} $4$-vector operator,
  \item[\textup{(PL2)}] each component $W_\mu$ commutes with all
    translations $P_\nu$,
  \item[\textup{(PL3)}] $\lim_{q\rightarrow 1} W_\mu = W_\mu^{q=1}$ as
    defined in~\eqref{eq:PauliLubanski1},
  \end{Proplist}
  is called a $q$-Pauli-Lubanski vector.
\end{Definition}
Obviously, (PL1)-(PL3) do not determine $W_\mu$ uniquely. For example,
we could multiply it by any $q$-polynomial which evaluates to 1 at
$q=1$. Property (PL1) tells us that the square of $W_\mu$ is a
$q$-Lorentz scalar, that is, commutes with all $h\in\Hcal$. As a
consequence of (PL2) this square commutes with all momenta. Therefore,
(PL1) and (PL2) together guarantee that the square of a
$q$-Pauli-Lubanski vector is a Casimir operator. The additional
property (PL3) is the obvious requirement that $W_\mu$ be a
$q$-deformation of the undeformed Pauli-Lubanski vector.

\subsection{Constructing the \textit{q}-Pauli-Lubanski Vector}

The only 4-vector operator we know so far is the 4-mo\-men\-tum
$P_\mu$. By construction, it is a \emph{left} 4-vector operator. Being
given a universal $\R$-matrix, there is a generic way to construct a
right tensor operator from of a given left tensor operator:
\begin{Proposition}
  \label{th:PL2}
  Let $T_j$ be a left $\rho$-tensor operator, that is, $\adL h \tr T_j
  = T_i \,\rho(h)^i{}_j$ for all $h\in\Hcal$, where $\rho$ is a finite
  representation of $\Hcal$. Let $\R$ be a universal $\R$-matrix of
  $\Hcal$. The set of operators
  \begin{equation}
  \label{eq:Sigma-def}
    \Sigma_\R(T_j) :=
    S^2(\R_{[1]})T_i\rho(\R_{[2]})^{i}{}_j 
  \end{equation}
  is a right $\rho$-tensor operator.
\end{Proposition}
\begin{proof}
  Abbreviating $\adL h \tr T_j \equiv h\tr T_j$, we have
\begin{equation}
\begin{split}
  \Sigma_\R(T_j)\tl\adR h
  &= S(h_{(1)}) S^2(\R_{[1]})(\R_{[2]}\tr T_j)h_{(2)}\\
  &= S(h_{(1)})S^2(\R_{[1]}) h_{(3)}
  \bigl( S^{-1}(h_{(2)})\R_{[2]} \tr T_j \bigr) \\
  &= S\bigl(S(\R_{[1]})h_{(1)}\bigr)h_{(3)}
  \bigl(S^{-1}(S(\R_{[2]})h_{(2)})\tr T_j \bigr) \\
  &= S\bigl(h_{(2)}S(\R_{[1]})\bigr)h_{(3)}
  \bigl(S^{-1}(h_{(1)}S(\R_{[2]}))\tr T_j \bigr)\\
  &= S^2(\R_{[1]}) S(h_{(2)})h_{(3)}
  \bigl(\R_{[2]}S^{-1}(h_{(1)})\tr T_j \bigr)\\
  &= S^2(\R_{[1]})(\R_{[2]}S^{-1}h\tr T_j)\\
  &= \Sigma_\R(T_i) \rho( S^{-1}h )^i{}_j \,.
\end{split}
\end{equation}
According to Eq.~\eqref{eq:TensorOp}, $\Sigma_\R(T_j)$ is indeed a
right $\rho$-tensor operator.
\end{proof}
This proposition tells us in particular, that $\Sigma_\R(P_\mu)$
satisfies (PL1). The next proposition takes care of (PL2).
\begin{Proposition}
  \label{th:PL3}
  Let $P_\mu$ be the momentum 4-vector, $\RI$ the antireal universal
  $\R$-matrix of the $q$-Lorentz algebra, and $\Sigma$ be defined as
  in Proposition~\ref{th:PL2}. Then
  \begin{equation}
    [\Sigma_\RI(P_\mu), P_\nu] =0 \,,\quad
    [\Sigma_{\R^{-1}_{\mathrm{I}\,21}}(P_\mu), P_\nu] =0
  \end{equation}
  for all $\mu$, $\nu$.
\end{Proposition}
\begin{proof}
  We denote the 4-vector representation by $\adL h \tr P_\nu = h\tr
  P_\nu = P_\mu \Lambda(h)^\mu{}_\nu$. Recall, that the commutation
  relations of the momenta can be written as
  \begin{equation}
    P_{\mu}P_{\nu} R_\mathrm{I}^{\nu\mu}{\!}_{\sigma\tau}
    = P_\sigma P_\tau \quad\Leftrightarrow\quad P_{\mu}P_{\nu}
    (R_\mathrm{I}^{-1})^{\mu\nu}{\!}_{\sigma\tau} = P_\tau P_\sigma \,,
  \end{equation}
  where the $R$-matrix $R_\mathrm{I}^{\mu\nu}{\!}_{\sigma\tau} =
  \Lambda(\R_{\mathrm{I}[1]})^{\mu}{}_\sigma
  \Lambda(\R_{\mathrm{I}[2]})^{\nu}{}_{\tau}$ is the 4-vector
  representation of the antireal universal $\R$-matrix $\RI$. Using
  the commutation relations between tensor operators and the Hopf
  algebra, we find
  \begin{equation}
  \begin{split}
    P_\sigma \,\Sigma_\RI(P_\tau) &= P_\sigma\,
    S^2(\R_{\mathrm{I}[1]})
    P_\nu \Lambda(\R_{\mathrm{I}[2]})^{\nu}{}_\tau \\
    &= S^2(\R_{\mathrm{I}[1']}) \bigl[ S(\R_{\mathrm{I}[1]})\tr
    P_\sigma \bigr] P_\nu \Lambda(\R_{\mathrm{I}[2]}
    \R_{\mathrm{I}[2']})^{\nu}{}_\tau \\
    &= S^2(\R_{\mathrm{I}[1']}) P_{\mu}
    \Lambda(\R^{-1}_{\mathrm{I}[1]})^{\mu}{}_{\sigma} P_\nu
    \Lambda(\R^{-1}_{\mathrm{I}[2]})^{\nu}{}_{\tau'}
    \Lambda(R_{\mathrm{I}[2']})^{\tau'}{}_\tau \\
    &= S^2(\R_{\mathrm{I}[1']}) \bigl[ P_{\mu} P_\nu
    (R_\mathrm{I}^{-1})^{\mu\nu}{\!}_{\sigma\tau'} \bigr]
    \Lambda(\R_{\mathrm{I}[2']})^{\tau'}{}_\tau \\
    &= S^2(\R_{\mathrm{I}[1']}) P_{\tau'} P_\sigma
    \Lambda(\R_{\mathrm{I}[2']})^{\tau'}{}_\tau \\
    &= \Sigma_\RI (P_\tau)\, P_\sigma \,.
  \end{split}
  \end{equation}
  On the second an third line we have used~\eqref{eq:Ph-rels},
  $\Delta(\R_{[1]})\otimes \R_{[2]} = \R_{[1]} \otimes \R_{[1']}
  \otimes \R_{[2]}\R_{[2']}$, and $S(\R_{[1]})\otimes \R_{[2]} =
  \R^{-1}_{[1]}\otimes \R^{-1}_{[2]}$.  The calculations for $\RI
  \rightarrow \R^{-1}_{\mathrm{I}\,21}$ are completely analogous.
\end{proof}
Propositions~\ref{th:PL2} and \ref{th:PL3} tell us, that both
$\Sigma_\RI(P_\mu)$ and $\Sigma_{\R^{-1}_{\mathrm{I}\,21}}(P_\mu)$
satisfy properties (PL1) and (PL2), respectively. In order to check
(PL3) we must find explicit expressions for $\Sigma_\RI(P_\mu)$ and
$\Sigma_{\R^{-1}_{\mathrm{I}\,21}}(P_\mu)$. This amounts to
calculating the $L$-matrices
\begin{equation}
  (L^\Lambda_{\mathrm{I}+})^\mu{}_\nu
  := \R_{\mathrm{I}[1]}\Lambda(\R_{\mathrm{I}[2]})^{\mu}{}_\nu
  \,,\quad (L^\Lambda_{\mathrm{I}-})^\mu{}_\nu
  := \R^{-1}_{\mathrm{I}[2]}
  \Lambda(\R^{-1}_{\mathrm{I}[1]})^{\mu}{}_\nu \,.
\end{equation}
For the 4-vector of these $L$-matrices we find
\begin{subequations}
\begin{align}
  (L^\Lambda_{\mathrm{I}+})^\mu{}_\nu &=
  \begin{pmatrix}
    1 & 0 & 0 & 0 \\
    0 & a^2 & b^2 & q^{\frac{1}{2}}[2]^{\frac{1}{2}} ab \\
    0 & c^2 & d^2 & q^{\frac{1}{2}}[2]^{\frac{1}{2}} cd \\
    0 & q^{\frac{1}{2}}[2]^{\frac{1}{2}} ac &
    q^{\frac{1}{2}}[2]^{\frac{1}{2}} bd & (1 + [2]bc)
  \end{pmatrix}\\
  (L^\Lambda_{\mathrm{I}-})^\mu{}_\nu &=
  \begin{pmatrix}
    W & \lambda K^{-1}J_- & \lambda K^{-1}J_+ & W-K^{-1} \\
    -q^{-1}\lambda J_+ & 1 & 0 & -q^{-1}\lambda J_+ \\
    -q \lambda J_- & 0 & 1 & -q\lambda J_-\\
    \lambda J_3 & -\lambda K^{-1}J_- & -\lambda K^{-1}J_+ & \lambda
    J_3 + K^{-1}
  \end{pmatrix}
\end{align}
\end{subequations}
with respect to the basis $\{0,-,+,3\}$. Observe, that
$(L^\Lambda_{\mathrm{I}+})^A{}_B$ is the 3-dimensional
corepresentation matrix of $\SUq^\op$. With a linear combination of
these two $L$-matrices we can satisfy (PL3).
\begin{Theorem}
\label{th:PL}
  The set of operators
  \begin{equation}
  \label{eq:PL-def}
  \begin{split}
    W_\nu &:=
    \lambda^{-1}\bigl[\Sigma_{\R^{-1}_{\mathrm{I}\,21}}(P_\nu)
    - \Sigma_{\RI}(P_\nu)\bigr] 
    = \lambda^{-1} S^2 \bigl[(L^\Lambda_{\mathrm{I}-})^\mu{}_\nu -
    (L^\Lambda_{\mathrm{I}+})^\mu{}_\nu \bigr] P_\mu
  \end{split}
  \end{equation}
  is a $q$-Pauli-Lubanski vector in the sense of
  Definition~\ref{th:PLdef}.
\end{Theorem}
\begin{proof}
  Properties (PL1) and (PL2) have been shown in
  Propositions~\ref{th:PL2} and \ref{th:PL3}, respectively. It remains
  to show (PL3). We note that the undeformed limit of the (left and
  right) Hopf adjoint action is the ordinary adjoint action. Hence,
  the limit $q\rightarrow 1$ preserves tensor operators. Since a
  4-vector is an \emph{irreducible} tensor operator it is sufficient
  to examine the limit of one component only. The limits of the other
  components follow by application of the adjoint action. We choose
  the zero component for which we have to show that
  \begin{equation}
    W_0 = \lambda^{-1} (W -1)P_0  + J_A P_B g^{AB}
    \stackrel{q \rightarrow 1}{\longrightarrow}
    W_0^{q=1} = J_A P_B g^{AB} \,.
  \end{equation}
  All there is to show is that
  \begin{equation}
  \begin{split}
    \lambda^{-1} (W -1) &= \lambda^{-1}\bigl( [2]^{-1}[q^{-1}K +
    qK^{-1}
    +\lambda^2 EF] -1 \bigr)\\
    &= \lambda^{-1}[2]^{-1}\bigl(q^{-1}K + qK^{-1} -[2]\bigr) +\lambda
    [2]^{-1} EF
  \end{split}
  \end{equation}
  vanishes for $q\rightarrow 1$. Clearly, the $\lambda [2]^{-1}EF$
  term of the last line vanishes. Using $K = q^H$ we get for the other
  terms
  \begin{multline}
    \lambda^{-1}(q^{-1}K + qK^{-1} -[2]) =\sum_{n=1}^{\infty}
    \frac{(q+(-1)^n q^{-1})
      (\ln q)^n}{\lambda n!}\, H^n\\
    = \ln q\, H + \frac{[2](\ln q)^2}{\lambda 2!} H^2 + \frac{(\ln
      q)^3}{3!} H^3 + \frac{[2](\ln q)^4}{\lambda 4!} H^4 + \ldots \,,
  \end{multline}
  which vanishes for $q\rightarrow 1$, since $\lim_{q\rightarrow
    1}\lambda^{-1} (\ln q)^n = 0$ for $n\geq 1$.
\end{proof}

\subsection{The Spin Casimir}

We proceed to calculate the spin Casimir
\begin{equation}
\label{eq:Casi1}
\begin{split}
  W^\tau W_\tau
  &= \eta^{\tau\nu} W_\nu W_\tau \\
  &= \lambda^{-2}\eta^{\tau\nu} S^2
  \bigl[(L^\Lambda_{\mathrm{I}-})^\mu{}_\nu -
  (L^\Lambda_{\mathrm{I}+})^\mu{}_\nu \bigr] P_\mu \,
  S^2\bigl[(L^\Lambda_{\mathrm{I}-})^\sigma{}_\tau -
  (L^\Lambda_{\mathrm{I}+})^\sigma{}_\tau \bigr] P_\sigma \\
  &= \lambda^{-2}\eta^{\tau\nu} S^2
  \bigl\{ \bigl[(L^\Lambda_{\mathrm{I}-})^\mu{}_\nu -
  (L^\Lambda_{\mathrm{I}+})^\mu{}_\nu \bigr] 
  \bigl[(L^\Lambda_{\mathrm{I}-})^\sigma{}_\tau -
  (L^\Lambda_{\mathrm{I}+})^\sigma{}_\tau \bigr] \bigr\}
  P_\sigma P_\mu \\
  &= \lambda^{-2}\eta^{\tau\nu} S^2 \bigl[
  (L^\Lambda_{\mathrm{I}-})^\mu{}_\nu
  (L^\Lambda_{\mathrm{I}-})^\sigma{}_\tau +
  (L^\Lambda_{\mathrm{I}+})^\mu{}_\nu
  (L^\Lambda_{\mathrm{I}+})^\sigma{}_\tau \\ &\qquad -
  (L^\Lambda_{\mathrm{I}-})^\mu{}_\nu
  (L^\Lambda_{\mathrm{I}+})^\sigma{}_\tau -
  (L^\Lambda_{\mathrm{I}+})^\mu{}_\nu
  (L^\Lambda_{\mathrm{I}-})^\sigma{}_\tau \bigr] 
  P_\sigma P_\mu \,.
\end{split}
\end{equation}
This can be further simplified. We first note that the
commutation relations of the $L$-matrices are such that
\begin{equation}
\begin{split}
  \eta^{\tau\nu}
  (L^\Lambda_{\mathrm{I}-})^\mu{}_\nu
  (L^\Lambda_{\mathrm{I}+})^\sigma{}_\tau P_\sigma P_\mu
  &= \eta^{\tau\nu} (R^{-1}_\mathrm{I})^{\tau'\nu'}{}_{\nu\tau}
  (L^\Lambda_{\mathrm{I}+})^{\mu'}{}_{\nu'}
  (L^\Lambda_{\mathrm{I}-})^{\sigma'}{}_{\tau'}
  R_\mathrm{I}^{\mu\sigma}{}_{\sigma'\mu'} P_\sigma P_\mu \\
  &= \eta^{\tau\nu}
  (L^\Lambda_{\mathrm{I}+})^\mu{}_\nu
  (L^\Lambda_{\mathrm{I}-})^\sigma{}_\tau P_\sigma P_\mu \,,
\end{split}
\end{equation}
where in the second step we have used the commutation
relations~\eqref{eq:MinkRel} of the momenta and that
$(R^{-1}_\mathrm{I})^{\tau'\nu'}{}_{\nu\tau} \eta^{\tau\nu} =
\eta^{\tau'\nu'}$. Moreover, using Eq.~\eqref{eq:metric1} one can see
that
\begin{equation}
 \eta^{\tau\nu} (L^\Lambda_{\mathrm{I}\pm})^\mu{}_\nu
  (L^\Lambda_{\mathrm{I}\pm})^\sigma{}_\tau
  = \eta^{\mu\sigma} \,. 
\end{equation}
With the last two results Eq.~\eqref{eq:Casi1} becomes
\begin{equation}
  W^\tau W_\tau = 2 \lambda^{-2} S^2 \bigl[
  \eta^{\mu\sigma} - \eta^{\tau\nu}
  (L^\Lambda_{\mathrm{I}+})^\mu{}_\nu
  (L^\Lambda_{\mathrm{I}-})^\sigma{}_\tau
  \bigr] P_\sigma P_\mu \,.
\end{equation}

\section{The Little Algebras}
\label{sec:qLittle}

\subsection{Little Algebras in the \textit{q}-Deformed Setting}
\label{sec:Little1}

In classical relativistic mechanics the state of motion of a free
particle is completely determined by its $4$-momentum. In quantum
mechanics particles can have an additional degree of freedom called
spin: Let us assume we have a free relativistic particle described by
an irreducible representation of the Poincar{\'e} algebra. We pick all
states with a given momentum,
\begin{equation}
  \Lcal_{\vec{p}} := \{ \Ket{\psi}\in \Lcal :
  P_\mu\Ket{\psi} = p_\mu\Ket{\psi} \} \,,
\end{equation}
where $\Lcal$ is the Hilbert space of the particle and
$\vec{p}=(p_\mu)$ is the $4$-vector of momentum eigenvalues. If the
state of the particle is \emph{not} uniquely determined by the
eigenvalues of the momentum, then the eigenspace $\Lcal_{\vec{p}}$
will be degenerate. In that case we need, besides the momentum
eigenvalues, an additional quantity to label the basis of our Hilbert
space uniquely. This additional degree of freedom is spin. The spin
symmetry is then the set of Lorentz transformations that leave the
momentum eigenvalues invariant and, hence, act on the spin degrees of
freedom only,
\begin{equation}
\label{eq:Little1}
  \mathcal{K}'_{\vec{p}} := \{ h\in\Hcal :
  P_\mu h \Ket{\psi} = p_\mu h \Ket{\psi}
  \text{ for all } \Ket{\psi}\in \Lcal_{\vec{p}} \} \,,
\end{equation}
where $\Hcal$ is the enveloping Lorentz algebra. In mathematical
terms, $\mathcal{K}'_{\vec{p}}$ is the stabilizer of $\Lcal_{\vec{p}}$.
Clearly, $\mathcal{K}'_{\vec{p}}$ is an algebra, called the little
algebra.

A priori, there are a lot of different little algebras for each
representation and each vector $p$ of momentum eigenvalues.  In the
undeformed case it turns out that for the physically relevant
representations (real mass) there are (up to isomorphism) only two
little algebras, depending on the mass being either positive or zero
\cite{Wigner:1939}. For positive mass we get the algebra of rotations,
$\mathcal{U}(\mathrm{su}_2)$, for zero mass an algebra which is
isomorphic to the algebra of rotations and translations of the
2-dimensional plane denoted by $\mathcal{U}(\mathrm{iso}_2)$.  The
proof that $\mathcal{K}'_{\vec{p}}$ does not depend on the particular
representation but on the mass does not generalize to the $q$-deformed
case: If we define for representations of the $q$-Poincar{\'e} algebra
the little algebra as in Eq.~\eqref{eq:Little1},
$\mathcal{K}'_{\vec{p}}$ for a spin-$\frac{1}{2}$ particle will not be
the same as for spin-1. We will therefore define the $q$-little
algebras differently.

In the undeformed case there is an alternative but equivalent
definition of the little algebras. $\mathcal{K}'_{\vec{p}}$ is the
algebra generated by the components of the $q$-Pauli-Lubanski vector
as defined in Eq.~\eqref{eq:PauliLubanski1} with the momentum
generators replaced by their eigenvalues. Let us formalize this to see
why this definition works and how it is generalized to the
$q$-deformed case.

Let $\chi_{\vec{p}}$ be the map that maps the momentum generators to the
eigenvalues, $\chi_{\vec{p}}(P_\mu) = p_\mu$. Being the restriction of a
representation, $\chi_{\vec{p}}$ must extend to a one dimensional
$*$-representation of the momentum algebra $\chi_{\vec{p}}: \Xcal
\rightarrow \mathbb{C}$, a non-trivial condition only in the
$q$-deformed case. Noting that every $a\in\Acal=\Xcal\rtimes\Hcal$ can
be uniquely written as $a = \sum_i h_i x_i$, where $h_i\in\Hcal$ and
$x_i\in\Xcal$, we can extend $\chi_{\vec{p}}$ to a linear map on all of
$\Acal$ by defining $\hat{\chi}_{\vec{p}}: \Acal \rightarrow \Hcal$ as
\begin{equation}
  \hat{\chi}_{\vec{p}}(\sum_i h_i x_i) := \sum_i h_i \chi_{\vec{p}}(x_i).
\end{equation}
The little algebra can now be alternatively defined as the unital
algebra generated by the images of the $q$-Pauli-Lubanski vector under
$\hat{\chi}_{\vec{p}}$,
\begin{equation}
\label{eq:Little2}
  \mathcal{K}_{\vec{p}}
  := \mathbb{C}\lrAngle{\hat{\chi}_{\vec{p}}(W_\mu)}\,.
\end{equation}
Why is this a reasonable definition? By construction the action of
every element of $\Acal$ on $\Lcal_{\vec{p}}$ is the same as of its
image under $\hat{\chi}_{\vec{p}}$. For any $\Ket{\psi}\in\Lcal_{\vec{p}}$
this means
\begin{equation}
  P_\mu \,\hat{\chi}_{\vec{p}}(W_\nu)\Ket{\psi}
  = \hat{\chi}_{\vec{p}}(P_\mu W_\nu)\Ket{\psi}
  = \hat{\chi}_{\vec{p}}(W_\nu P_\mu)\Ket{\psi}
  = p_\mu \,\hat{\chi}_{\vec{p}}(W_\nu)\Ket{\psi} \,,
\end{equation}
which shows that $\mathcal{K}_{\vec{p}} \subset
\mathcal{K}'_{\vec{p}}$. It still could happen, that
$\mathcal{K}_{\vec{p}}$ is strictly smaller than
$\mathcal{K}'_{\vec{p}}$.  In the undeformed case there are theorems
\cite{Blattner:1969,Dixmier} telling us that this cannot happen, so we
really have $\mathcal{K}_{\vec{p}} = \mathcal{K}'_{\vec{p}}$. For the
$q$-deformed case no such theorem is known \cite{Schneider:1994}.
However, if there were more generators in the stabilizer of some
momentum eigenspace they would have to vanish for $q\rightarrow 1$. In
this sense Eq.~\eqref{eq:Little2} with the $q$-deformed Pauli-Lubanski
vector can be considered to define the $q$-deformed little algebras.

\subsection{Computation of the \textit{q}-Little Algebras}
\label{sec:Little2}

To begin the explicit calculation of the $q$-deformed little algebras,
we need to figure out if there are eigenstates of $q$-momentum at all.
That is, we want to determine the one-dimensional $*$-representations
of $\Xcal = \Mink$, that is the homomorphisms of $*$-Algebras
$\chi:\Xcal \mapsto \mathbb{C}$. Let us again denote the eigenvalues
of the generators by lower case letters $p_\mu := \chi(P_\mu)$.
According to Eq.~\eqref{eq:Pstar}, we must have $p_0$, $p_3$ real and
$p_+^* = -q p_-$ for $\chi$ to be a $*$-map. To find the conditions
for $\chi$ to be a homomorphism of algebras, we apply $\chi$ to the
relations~\eqref{eq:PP-Rel} of $\Xcal = \Mink$, yielding $p_A(p_0-p_3)
= 0$. There are two cases. The first is $p_0 \neq p_3$, which
immediately leads to $p_A = 0$, and $p_0 = \pm m$. The second case is
$p_0 = p_3$, leading to $m^2 = - \lvert p_- \rvert^2 - \lvert p_+
\rvert^2$, where, if the mass $m$ is to be real, we must have $p_\pm =
0$.

In summary, for real mass $m$ we have a massive and a massless type of
momentum eigenstate with eigenvalues given by
\begin{equation}
  (p_0,p_-,p_+,p_3) =
  \begin{cases}
    (\pm m ,0,0,0), & m > 0 \\
    (k ,0,0,k), & m = 0,\, k\in\mathbb{R}
  \end{cases}
\end{equation}
According to Eq.~\eqref{eq:Little2} we now have to replace the momenta
in the definition~\eqref{eq:PL-def} of the $q$-Pauli-Lubanski vector
with these eigenvalues.

For \textbf{the massive case} we get
\begin{equation}
\begin{aligned}
  \hat{\chi}_{\vec{p}}(W_0) &= \lambda^{-1} (W -1)m \\
  \hat{\chi}_{\vec{p}}(W_-) &= J_- K^{-1} m \\
  \hat{\chi}_{\vec{p}}(W_+) &= J_+ K^{-1} m \\
  \hat{\chi}_{\vec{p}}(W_3) &= \lambda^{-1} (W - K^{-1}) m \,,
\end{aligned}
\end{equation}
so the set of generators of the little algebra is essentially
$\{W,K^{-1},J_\pm K^{-1}\}$. Since $K^{-1}$ stabilizes the momentum
eigenspace, so does its inverse $K$. Hence, it is safe to add $K$ to
the little algebra which would exist, anyway, as operator within a
representation. We thus get
\begin{equation}
  \mathcal{K}_{(m ,0,0,0)} = \suq \,,
\end{equation}
completely analogous to the undeformed case.

\textbf{The massless case} is more interesting. Replacing the momentum
generators with $(P_0,P_-,P_+,P_3) \rightarrow (k ,0,0,k)$ we get
\begin{equation}
\begin{aligned}
  \hat{\chi}_{\vec{p}}(W_0) &= \lambda^{-1} (K -1) k \\
  \hat{\chi}_{\vec{p}}(W_-) &=
  -\lambda^{-1} q^{-\frac{3}{2}}[2]^{\frac{1}{2}} ac \, k\\
  \hat{\chi}_{\vec{p}}(W_+) &=
  -\lambda^{-1}  q^{\frac{5}{2}}[2]^{\frac{1}{2}} bd \, k\\
  \hat{\chi}_{\vec{p}}(W_3) &= \lambda^{-1}\bigl( K - (1 +[2] bc) \bigr)
  k \,.
\end{aligned}
\end{equation}
The set of generators of this little algebra is essentially $\{K, ac,
bd, bc \}$. The commutation relations of these generators can be
written more conveniently in terms of $K$ and $N_A :=
(L^\Lambda_{\mathrm{I}+})^3{}_A$, that is 
\begin{xalignat}{3}
  N_- &= q^{\frac{1}{2}}[2]^{\frac{1}{2}} ac \,,& N_+ &=
  q^{\frac{1}{2}}[2]^{\frac{1}{2}} bd \,,& N_3 &= 1 +[2] bc \,.
\end{xalignat} 
The commutation relations are
\begin{xalignat}{3}
  N_B N_A\, \varepsilon^{AB}{}_C &= -\lambda N_C \,,& N_A N_B\, g^{BA}
  &= 1 \,,& K N_A &= q^{-2A} N_A K \,,
\end{xalignat}
with conjugation $N_A^* = N_B\, g^{BA}$, $K^* = K$. In words: The
$N_A$ generate the opposite algebra of a unit quantum sphere,
$\qsphere$ \cite{Podles:1987}. $K$, the generator of $\uq$, acts on
$N_A$ as on a right $3$-vector operator. In total we have
\begin{equation}
  \mathcal{K}_{(k ,0,0,k)} = \uq \ltimes \qsphere \,.
\end{equation}
As opposed to the massive case, this is no Hopf algebra. However,
since $L$-matrices are multiplicative, that is,
$\Delta[(L^\Lambda_{\mathrm{I}+})^\mu{}_\sigma] =
(L^\Lambda_{\mathrm{I}+})^\mu{}_\nu \otimes
(L^\Lambda_{\mathrm{I}+})^\nu{}_\sigma$, we have
\begin{equation}
 \Delta(N_{B}) = N_A \otimes
 (L^\Lambda_{\mathrm{I}+})^A{}_B\,,
\end{equation}
hence, $\mathcal{K}_{(k,0,0,k)}$ is a right coideal.

\appendix

\section{Appendix: The \textit{q}-Poincar\'e Algebra}
\label{sec:AppPoin}

The Hopf-$*$ algebra generated by $E$, $F$, $K$, and $K^{-1}$ with
relations
\begin{equation}
\begin{gathered}
  KK^{-1}= 1 = K^{-1}K \,,\quad KEK^{-1} = q^2 E \,,\\
  KFK^{-1} = q^{-2}F\,,\quad [E,F]= \lambda^{-1}(K-K^{-1})\,,
\end{gathered}
\end{equation}
Hopf structure
\begin{equation}
\begin{gathered}
  \Delta(E) = E\otimes K + 1\otimes E \,,\quad
  \Delta(F) = F\otimes 1 + K^{-1}\otimes F \,, \\
  \Delta(K) = K \otimes K \,,\quad
  \varepsilon(E)=0=\varepsilon(F) \,,\quad \varepsilon(K)=1\,,\\
  S(E) = -EK^{-1} \,,\quad S(F) = -KF \,,\quad S(K)=K^{-1}\,,
\end{gathered}
\end{equation}
and $*$-structure
\begin{equation}
  E^* = FK \,, F^* = K^{-1}E \,, K^* = K
\end{equation}
is called $\suq$, the $q$-deformation of the enveloping algebra $\su$.
  
The set of generators $\{J_A\} = \{J_-,J_3,J_+\}$ of $\suq$ defined as
\begin{equation}
\label{eq:DefJ}
\begin{aligned}
  J_{-} &:= q[2]^{-\frac{1}{2}}KF \\
  J_3   &:= [2]^{-1} (q^{-1}EF-qFE) \\
  J_{+} &:= -[2]^{-\frac{1}{2}}E
\end{aligned}
\end{equation}
is the left 3-vector operator of angular momentum. The center of
$\suq$ is generated by
\begin{equation}
\label{eq:DefW}
  W := K - \lambda J_3 = K - \lambda [2]^{-1} (q^{-1}EF-qFE)\,,
\end{equation}
the Casimir operator of angular momentum. $W$ is related to $J_A$ by
\begin{equation}
  W^2 -1 = \lambda^2(J_3^2 - q^{-1} J_- J_+ -  q J_+ J_-)
  = \lambda^2J_A J_B g^{AB} \,,
\end{equation}
thus defining the 3-metric $g^{AB}$, by which we raise 3-vector
indices $J^A = g^{AB}J_B$. It is also useful to define an
$\varepsilon$-tensor
\begin{xalignat}{3}
  \varepsilon^{-3}{}_- &= q^{-1} & \varepsilon^{3-}{}_- &= -q && \notag\\
  \varepsilon^{-+}{}_3 &= 1 & \varepsilon^{+-}{}_3 &= -1 &
  \varepsilon^{33}{}_3 &= -\lambda \\
  \varepsilon^{3+}{}_+ &= q^{-1} & \varepsilon^{+3}{}_+ &= -q \,,\notag &&
\end{xalignat}

The Hopf-$*$ algebra generated by the $2\times 2$-matrix of generators
$B^i{}_j = (\begin{smallmatrix}a&b\\c&d\end{smallmatrix})$ with
relations
\begin{equation}
\label{eq:SUopRel}
\begin{gathered}
  ba=qab, \quad ca=qac, \quad db = q bd, \quad dc = q cd \\
  bc = cb, \quad da - ad = (q - q^{-1}) bc, \quad da-qbc=1 \,,
\end{gathered}
\end{equation}
coproduct $\Delta (B^{i}{}_{k}) = B^{i}{}_{j}\otimes B^{j}{}_{k}$
(summation over j), counit $\varepsilon (B^{i}{}_{j}) =
\delta^{i}_{j}$, antipode and $*$-structure
\begin{xalignat}{2}
  S\begin{pmatrix} a & b \\ c & d \end{pmatrix} &=
  \begin{pmatrix} d & -qb \\ -q^{-1}c & a \end{pmatrix}, &
  \begin{pmatrix} a & b \\ c & d \end{pmatrix}^* &=
  \begin{pmatrix} d & -q^{-1}c \\ -qb & a \end{pmatrix},
\end{xalignat}
is $\SUq^\op$, the opposite algebra of the quantum group $\SUq$.

The Hopf-$*$ algebra generated by the Hopf-$*$ sub-algebras $\suq$ and
$\SUq^\op$ with cross commutation relations
\begin{equation}
\begin{aligned}
  \begin{pmatrix} a & b \\ c & d \end{pmatrix} E
  &= \begin{pmatrix} q E a - q^{\frac{3}{2}} b & q^{-1}Eb \\
    qEc+q^{\frac{3}{2}}Ka- q^{\frac{3}{2}}d &
    q^{-1}Ed+q^{-\frac{1}{2}}Kb \end{pmatrix} \\
  \begin{pmatrix} a & b \\ c & d \end{pmatrix} F
  &= \begin{pmatrix} q F a + q^{-\frac{1}{2}}c &
    qFb-q^{-\frac{1}{2}}K^{-1}a + q^{-\frac{1}{2}}d \\
    q^{-1}Fc & q^{-1}Fd-q^{-\frac{5}{2}}K^{-1}c \end{pmatrix} \\
  \begin{pmatrix} a & b \\ c & d \end{pmatrix} K
  &= K \begin{pmatrix} a & q^{-2}b \\ q^{2} c & d\end{pmatrix},
\end{aligned}
\end{equation}
which is the Drinfeld double of $\suq$ and $\SUq^\op$, is the
$q$-Lorentz algebra $\Hcal = \slC$ \cite{Podles:1990}.

The $*$-algebra generated by $P_0$, $P_-$, $P_+$, $P_3$ with
commutation relations
\begin{xalignat}{2}
\label{eq:PP-Rel}
  P_0 P_A &= P_A P_0 \,,&
  P_A P_B\,\varepsilon^{AB}{}_{C} &= -\lambda P_0 P_C 
\end{xalignat}
and $*$-structure 
\begin{equation}
\label{eq:Pstar}
  P_0^* = P_0,\quad P_-^* = -q^{-1} P_+,\quad
  P_+^* = -q,\quad P_-,\quad P_3^* = P_3
\end{equation}
is the $q$-Minkowski space algebra $\Xcal = \Mink$. The center of
$\Mink$ is generated by the mass Casimir
\begin{equation}
\label{eq:FourMetric}
  m^2 := P_\mu P_\nu \eta^{\mu\nu}
  = P_0^2  + q^{-1} P_- P_+ +  q P_+ P_- - P_3^2 \,,
\end{equation}
thus defining the 4-metric $\eta^{\mu\nu}$. It is related to the
3-metric by $\eta^{AB} = - g^{AB}$ for $A,B\in\{-,+,3\}$.

The commutation relations of $\Xcal = \Mink$ are consistent with the
4-vector action $h\tr P_\nu = P_\mu \Lambda(h)^\nu{}_\mu$ of $\Hcal$
on $\Xcal$. $\Lambda$ is defined on the generators of rotations as
\begin{equation}
  \Lambda(J_A) =
  \begin{pmatrix} \rho^0(J_A) & 0 \\ 0 & \rho^1(J_A) \end{pmatrix} =
  \begin{pmatrix} 0 & 0 \\ 0 & \varepsilon_A{}^B{}_C \end{pmatrix} \,,
\end{equation}
where $\rho^0$ and $\rho^1$ are the spin-0 and the spin-1
representations of $\suq$, respectively. On the boost generators
$\Lambda$ is given by
\begin{subequations}
\begin{align}
  \Lambda(a) &=
  \begin{pmatrix}
    [4][2]^{-2} & 0 & 0 & q\lambda [2]^{-1} \\
    0 & 1 & 0 & 0 \\
    0 & 0 & 1 & 0 \\
    q^{-1}\lambda [2]^{-1} & 0 & 0 & 2 [2]^{-1} \\
  \end{pmatrix} \displaybreak[0] \\ 
  \Lambda(b) &= q^{-\frac{1}{2}}\lambda[2]^{-\frac{1}{2}}
  \begin{pmatrix}
    0 & -1 & 0 & 0 \\
    0 & 0 & 0 & 0 \\
    1 & 0 & 0 & 1 \\
    0 & 1 & 0 & 0 \\
  \end{pmatrix} \displaybreak[0] \\
  \Lambda(c) &= -q^{\frac{1}{2}}\lambda[2]^{-\frac{1}{2}}
  \begin{pmatrix}
    0 & 0 & -1 & 0 \\
    1 & 0 & 0 & 1 \\
    0 & 0 & 0 & 0 \\
    0 & 0 & 1 & 0 \\
  \end{pmatrix} \displaybreak[0] \\
  \Lambda(d) &=
  \begin{pmatrix}
    [4][2]^{-2} & 0 & 0 & -q^{-1}\lambda [2]^{-1} \\
    0 & 1 & 0 & 0 \\
    0 & 0 & 1 & 0 \\
    -q \lambda [2]^{-1} & 0 & 0 & 2 [2]^{-1}
  \end{pmatrix}
\end{align}
\end{subequations}
with respect to the $\{0,-,+,3\}$ basis. It has the property
\begin{equation}
\label{eq:metric1}
  \eta^{\nu\nu'} \Lambda(h)^{\mu'}{}_{\nu'} \eta_{\mu'\mu} =
  \Lambda(Sh)^{\nu}{}_{\mu}
\end{equation}
for all $h\in\Hcal$.

Finally, the $q$-Poincar\'e algebra is the $*$-algebra generated by
the $q$-Lorentz algebra $\Hcal = \slC$ and the $q$-Minkowski algebra
$\Xcal = \Mink$ with cross commutation relations
\begin{equation}
\label{eq:Ph-rels}
  hP_\nu = P_\mu \Lambda(h_{(1)})^\mu{}_\nu \, h_{(2)}
  \quad\Leftrightarrow\quad
  P_\nu\, h = h_{(2)} P_\mu
  \Lambda\bigl(S^{-1}(h_{(1)}) \bigr)^\mu{}_\nu \,.
\end{equation}
More details and mathematical background information has been
compiled in \cite{Blohmann}.

\subsection*{Acknowledgement}
This work was supported by the Studienstiftung des deut\-schen Volkes.

\end{document}